\documentclass{naw-latex}
               {\endgraf\smallbreak
                \noindent\textbf{Definitie}\sffamily}%
               {\endgraf\smallskip}

\bibitem{MR512106}{F.~Ledrappier,  `Un champ markovien peut \^etre d'entropie nulle et
  m\'elangeant', \emph{C. R. Acad. Sci. Paris S\'er. A-B} \textbf{287} (1978),
  no.~7, A561--A563.}

\bibitem{MR1048237}{H. P. Schlickewei, `An explicit upper bound for the number of solutions of the {$S$}-unit equation',
\emph{J. Reine Angew. Math.}, \textbf{406} (1990),
109--120.}

\bibitem{MR990217}{K. Schmidt, `Mixing automorphisms of compact groups and a theorem by Kurt Mahler',
\emph{Pacific J. Math.} \textbf{137} (1989),
no~2, 371--385.}

\bibitem{MR97c:28041}{K.~Schmidt, \emph{Dynamical systems of algebraic origin} (Birkh\"auser Verlag,
  Basel, 1995).}

\bibitem{MR1193598}{K.~Schmidt and T.~Ward,
`Mixing automorphisms of compact groups and a theorem
of {S}chlickewei', \emph{Invent. Math.} \textbf{111} (1993), no.~1, 69--76.}

\begin{document}

\StartArtikel[Titel={Mixing, wine, and serendipity},
          AuteurA={T. Ward},
          AdresA={Ziff Building Floor 13\crlf
                  University of Leeds\crlf
                  Leeds\crlf
                  LS2 9JT, UK},
          EmailA={t.b.ward@leeds.ac.uk},
          kolommen={3}%%%%{3} mag ook
	  ]
	
\StartLeadIn

\StopLeadIn

During the~$1980$s, while I was a graduate student at Warwick University
under the supervision of Klaus Schmidt, a specific kind of algebraic
dynamical system was emerging as a surprisingly rich and relatively
unexplored field. In hindsight, a small shift in how a key example
constructed by Ledrappier
is thought of might have predicted some of this, but prediction with the
benefit of hindsight is a little too easy.

Mixing is a mathematical version of the idea of, well, mixing.
If two ingredients of a cocktail are poured carefully into a glass --- so carefully
that perhaps they form individual layers --- then the action of a
stirrer is `mixing' if after some time every mouthful tastes the same. That is,
every part of the glass has the ingredients in the same proportion up to a
negligible error. This becomes a mathematical concept by noticing that
the volume may be viewed as a measure on the space
consisting of the contents of the glass, and the action of the stirrer might
be thought of as iteration of a map that preserves that measure (unless it is
being stirred using a straw, and the person stirring is taking a crafty
sip every now and then). Avoiding all the interesting and subtle physical
and chemical issues involved --- particularly egregious in the circumstances --- we
might as well assume the action of stirring is invertible, and for
mathematicians the resulting structure of a measure-preserving
action of the integers might as well be an action of any group.
Having no wish to trip up on any measure theory,
let's say that the group is countable. So here is mixing:
if a countable group~$G$ acts by transformations
preserving a measure~$\mu$
on a probability space, then it is
called mixing if~$\mu(A\cap gB)$ converges to~$\mu(A)\mu(B)$
as~$g$ `goes to infinity' in~$G$. And why not be ambitious? Mixing
on~$(k+1)$ sets (or mixing of order~$k$) means that for any measurable sets~$A_0,\dots,A_k$
the measure
\[
\mu(A_0\cap g_1A_1\cap\cdots\cap g_kA_k)
\]
of the intersection converges to~$\prod_{j=0}^{k}\mu(A_j)$
as the group elements~$g_j$ go to infinity and move apart from each other.
So here is a mathematical question: given a measure-preserving
action of a countable group, determine if it is mixing on~$k$ sets for a given~$k$.
When~$G=\mathbb{Z}$ it is a long-standing question of Rokhlin as
to whether mixing on~$2$ sets forces mixing on~$3$ sets.

Which brings us to Ledrappier's example~\cite{MR512106}
(simplified for convenience from his harmonic condition example):
let~$X$ be the subset of~$\{0,1\}^{\mathbb{Z}^2}$ consisting
of the points~$x$ satisfying~$x_{s,t}=x_{s+1,t}+x_{s,t+1}$ modulo~$2$
for every~$(s,t)\in\mathbb{Z}^2$.
This is a compact group, and the shift in~$\mathbb{Z}^2$ defines
an action of~$\mathbb{Z}^2$ that preserves the natural Haar measure.
The system is easily shown to be mixing on~$2$ sets,
but the fact that the relation~$x_{s,t}=x_{s+2^n,t}+x_{s,t+2^n}$ modulo~$2$
holds for all~$n\ge1$ (a direct consequence of the
properties of the Frobenius~$a\mapsto a^2$ modulo~$2$ under iteration)
forces a correlation between triples of sets
separated by arbitrarily large distances --- failure of mixing on~$3$ sets.
Ledrappier also pointed out that any system like this built from
automorphisms of compact groups has a property called `Lebesgue spectrum'.
A productive shift in perspective is to think of this system as the dual
group of the module~$\mathbb{Z}[u_1^{\pm1},u_2^{\pm1}]/\langle 1+u_1+u_2\rangle$, with
the action of~$(a,b)\in\mathbb{Z}^2$ dual to multiplication by~$u_1^au_2^b$.
Thus a version of the mixing question becomes this: describe the
mixing properties of such a system built from a module~$M$ over the ring~$R=\mathbb{Z}[u_1^{\pm1},\dots,u_d^{\pm1}]$
in terms of properties of the module~$M$ --- in the certain knowledge
that the answer is non-trivial because it is for Ledrappier's example.

Work of Kitchens and Schmidt~\cite{MR97c:28041} probed the
mixing properties of systems whose compact group is
zero-dimensional, uncovering a complex collection of
properties leading to many interesting questions.
Schmidt~\cite{MR990217} also showed that
the way in which a `shape' produced by the Frobenius
automorphism witnesses failure of
higher-order mixing as seen by Ledrappier could not
take place if the compact group~$X$ is connected. That is,
for a mixing~$\mathbb{Z}^d$ action by automorphisms of
a compact connected group, choosing the times~$g_1,\dots,g_k$
to be dilates of a fixed shape in~$\mathbb{Z}^d$ would never
show failure of mixing, raising the question: for
these connected systems, does mixing imply mixing of all orders?

By~$1991$ I was working at Ohio State University, and
we were notified that some
duplicate journals were being discarded. As life was then full
of time for mathematics, I went
into the basement and leafed through piles of journals
in recycling bins, tearing out
any articles that looked vaguely interesting.
I piled these up, and left
to attend a workshop at
CIRM in Luminy.

In that beautiful place, Klaus Schmidt reminded
Doug Lind and me of this open problem over a splendid meal. Perhaps with
the assistance of the generous provision of wine, I became
sure that I had an argument, essentially using the Lebesque
spectrum property, that proved mixing of all orders for these
connected systems. Not for the first, and not for the
last, time, Klaus let me whitter on for some time as we walked
under the pine trees before
politely pointing out that my suggested argument applied
unchanged to Ledrappier's example.

Flying back to Columbus, the problem was
firmly in my mind. The~`$\times2,\times3$' system,
itself studied for other reasons, was
the natural start.
Ledrappier's salutary example showed that the result sought really
couldn't come from the familiar toolbox of spectral
or entropy methods. Via Fourier analysis of indicator
functions of sets,
it seemed to come
down to this: what can you say about
solutions of~$\sum_{j=1}^{k}a_jx_j=1$
in a number field, where the variables~$x_j$
come from a finitely-generated multiplicative subgroup?
For~`$\times2,\times3$' the field would be~$\mathbb{Q}$,
and the multiplicative subgroup~$\{2^a3^b\mid a,b\in\mathbb{Z}\}$.
Failure of mixing of all orders seemed roughly equivalent
to equations of this shape having too many --- infinitely many ---
non-trivially different solutions. Trivially
different solutions abound if a sub-sum vanishes,
because that vanishing sub-sum can be scaled by
powers of~$2$ and~$3$ arbitrarily.

After a few days back in Columbus, I sorted through the pile of
torn-out papers on my desk. One was
a (then) recent paper of Schlickewei~\cite{MR1048237}
with a form of `$S$-unit theorem'. For the
finite-dimensional case
at hand (it turned out later that the topological
dimension of the compact group~$X$ plays a role)
a simple reduction argument was suddenly completely clear.
If a~$\mathbb{Z}^d$-action by automorphisms of
a compact connected abelian group fails to
be mixing on~$k\ge3$ sets, then by the Fourier analysis
argument there is a linear equation with~$k$
terms over a field of characteristic zero that has infinitely many
distinct solutions lying in a multiplicative group with~$d$
generators (corresponding to the automorphisms defining the action).
By the~$S$-unit theorem, this is impossible unless
infinitely many of them come from a vanishing sub-sum: a linear
equation with~$j<k$ terms. But
finding infinitely many solutions for that shorter linear
equation is a witness to failure of mixing on~$j<k$ sets.
Thus mixing on~$2$ sets implies mixing of all orders.
Some cleaning up was needed --- algebra to reduce to
cyclic modules, and a more subtle process needed to deal with
infinite-dimensional compact groups which do not
readily permit the translation into statements in number
fields --- but this quickly led to the proof
of the full case with Klaus Schmidt~\cite{MR1193598}.

There are some lessons to take from this strange
coincidence and happy resolution. Certainly ideas produced under
the influence of wine may
eventually face the sobering reality
of counter-examples --- but can be motivating
nonetheless. More importantly,
Tramezzino's tale
{\it Peregrinaggio di tre giovani figliuoli del re di Serendippo}
in which ``accidents and sagacity'' play such a role
still has something for us.
The rush of modern academic life, the growing use of online
journals and their sophisticated and well-intentioned nudging
towards related articles, the demarcation of subject areas,
the overwhelming growth in the volume
of the mathematical literature --- all make the benefits
of serendipity less easy to access.
If you have the good fortune of time to spend on
mathematics, spend some of it on the {not} `suggested
article', on the articles that readers of your
article are {not} `also reading', and on articles with
the {wrong} subject classification codes.

\completepublications

\end{document}